\documentclass[11pt]{article} 
 
\usepackage{amsmath} 
\usepackage{amssymb} 
\usepackage{latexsym}

\usepackage{latexsym} 
\usepackage{amsfonts} 

\usepackage{amsmath,amsthm,amssymb}

\newtheorem{theor}{Theorem} 
 
\newtheorem{corol}{Corollary}
\newtheorem{thm}{Theorem}[section] 
\newtheorem{lem}[thm]{Lemma} 
\newtheorem{cor}[thm]{Corollary} 
 
\newtheorem{prop}[thm]{Proposition} 
\theoremstyle{definition} 
\newtheorem{defn}[thm]{Definition} 
 
\newtheorem{conv}[thm]{Convention} 
\newtheorem{rem}[thm]{Remark}

\begin{document}

\title{Clusters, currents and Whitehead's algorithm} 
 
\author{\Large Ilya Kapovich}

 
\date{} 
 
\maketitle

\begin{abstract} 
 
Using geodesic currents, we provide a theoretical justification 
for some of the experimental results, obtained by Haralick, 
Miasnikov and Myasnikov via pattern recognition methods, regarding 
the behavior of Whitehead's algorithm on non-minimal inputs. In 
particular we prove that the images of ``random" elements of a 
free group $F$ under the automorphisms of $F$ form ``clusters" 
that share similar normalized Whitehead graphs and similar 
behavior with respect to Whitehead's algorithm. 
\end{abstract}

 


\section{Introduction}\label{intro}

The \emph{automorphism problem} for a free group $F=F(a_1,\dots, 
a_k)$, where $k\ge 2$, asks, given two arbitrary elements $u,v\in 
F$, whether there exists $\phi\in Aut(F)$ such that $\phi(u)=v$. 
In a classic 1936 paper~\cite{Wh} Whitehead provided an algorithm 
solving the automorphism problem. He introduced a special finite 
generating set of $Aut(F)$, consisting of the so-called 
\emph{Whitehead automorphisms}. He proved that if $u\in F$ is a 
cyclically reduced word that is not shortest in its 
$Aut(F)$-orbit, then there exists a Whitehead automorphism $\tau$ 
such that $\tau(u)$ has smaller cyclically reduced length than 
$\tau$. This provides a quadratic time algorithm for finding a 
\emph{minimal} element in the orbit $Aut(F)f$ for any $f\in F$, 
that is, the element of smallest length in $Aut(F)f$. Namely, 
first cyclically reduce $f$ to get $f'\in F$, and then check if 
there is a Whitehead automorphism $\tau$ that decreases the 
cyclically reduced length of $f'$. If not, then $f'$ is minimal. 
If yes, replace $f'$ by $\tau(f)$ and then repeat the entire step. 
Whitehead also proved that if $u,v\in F$ are cyclically reduced 
minimal elements of the same length, then $v\in Aut(F)u$ if and 
only if there exists a chain of Whitehead automorphisms taking $u$ 
to $v$ and such that the cyclically reduced length is constant 
throughout the chain. Together with the above procedure for 
computing minimal representatives, this provides an algorithm for 
solving the automorphism problem that runs in at most exponential 
time in terms of $|u|+|v|$. The second, ``hard'' part of 
Whitehead's algorithm, has an a priori exponential time upper 
bound for the running time, although in practice the algorithm 
appears to always terminate much faster. 
 
Since this 1936 paper of Whitehead there has been a great deal of 
work on the study of the automorphism problem and of Whitehead's 
algorithm (e.g. see the recent paper of Lee~\cite{Lee}). However, 
even now, 70 years later, it is still not known what the precise 
complexity of Whitehead's algorithm is or if there exists a 
polynomial time algorithm for solving the automorphism problem in 
a free group. The only well-understood case is $k=2$, where it is 
known that the automorphism problem is indeed solvable in 
polynomial time~\cite{MS,Khan}. 
 
A recent paper of Kapovich, Schupp and Shpilrain~\cite{KSS} proves 
that for any $k\ge 2$ Whitehead's algorithm has linear time 
generic-case complexity. It turns out that ``random'' cyclically 
reduced elements of $F$ are already minimal, so that the first 
(minimization) part of Whitehead's algorithm terminates in a 
single step. Moreover, even the second ``hard'' part of the 
algorithm is also proved in \cite{KSS} to run in at most linear 
time in this case. 
 
It is therefore interesting to understand the behavior of Whitehead's 
algorithm on non-minimal inputs that are also generated via some 
natural probabilistic process. A.~D.~Miasnikov, A.~G.~Myasnikov and 
R.~Haralick~\cite{HMM1,HMM2,HMM3}, via pattern recognition methods, 
experimentally discovered some interesting features of the behavior of 
Whitehead's algorithm in this set-up. Before discussing their 
observations, we need to fix some notations.

\begin{conv} 
  For the remainder of the paper let $F=F(A)$ be a free group with a 
  fixed free basis $A=\{a_1,\dots, a_k\}$, where $k\ge 2$. Let 
  $X=\Gamma(F,A)$ be the Cayley graph of $F$ with respect to $A$, so 
  that $X$ is a $(2k)$-regular tree. 
 
  Denote $\Sigma=A\cup A^{-1}=\{a_1,\dots, a_k, a_1^{-1}, \dots, 
  a_k^{-1}\}$. 
 
  For a word $w$ in $\Sigma^*$ we will denote the length of $w$ by 
  $|w|$. A word $w\in \Sigma^*$ is said to be \emph{reduced} if $w$ is 
  freely reduced in $F$, that is $w$ does not contain subwords of the 
  form $a_ia_i^{-1}$ or $a_i^{-1}a_i$. A word $w$ is \emph{cyclically 
    reduced} if all cyclic permutations of $w$ are reduced. (In 
  particular $w$ itself is reduced.) We denote by $C$ the set of all 
  nontrivial cyclically reduced words in $F$. 
 
  Since every element of $F$ can be uniquely represented by a freely 
  reduced word, we identify elements of $F$ and freely reduced words. 
  Any freely reduced element $w$ can be uniquely decomposed as a 
  concatenation $w = vuv^{-1}$ where $u$ is a cyclically reduced word. 
  The word $u$ is called the \emph{cyclically reduced form of $w$} and 
  $||w||:=|u|$ is the \emph{cyclic length} of $w$. 
\end{conv} 
 
Some of the experimental conclusions of A.~D.~Miasnikov, 
A.~G.~Myasnikov and R.~Haralick, described in detail in \cite{HMM3}, 
can be summarized as follows.  First take a large sample of long 
random cyclically reduced words $W_1$ in $F$.  If there are any 
non-minimal elements, apply Whitehead's algorithm and replace them by 
their minimal representatives. The resulting set $W_2$ consists of 
only minimal words. By the results of~\cite{KSS} most of elements of 
$W_1$ are already minimal and therefore the difference between $W_1$ 
and $W_2$ will be very small and can be disregarded. 
 
Then some of the elements $w$ of $W_2$ (again usually chosen at 
random) are replaced by $\phi_w(w)$ where $\phi_w$ comes from some 
finite collection $\Phi$ of automorphisms chosen so that 
$||w||<||\phi_w(w)||$. The resulting set $W_3$ thus contains both 
minimal and non-minimal elements. Some of the observed results were 
that: 
 
\begin{itemize} 
\item The non-minimal elements of the set $W_3$ formed several 
  ``clusters". 
 
\item For each ``cluster" $\mathcal C$ all the elements of 
$\mathcal C$ 
  had approximately the same normalized Whitehead graphs. 
 
\item Moreover, for each ``cluster" $\mathcal C$ there was a 
Whitehead 
  automorphism $\tau$ such that for all $w\in \mathcal C$ 
$$ 
||\tau(w)||<||w||. 
$$ 
(In fact, often, depending on how $\Phi$ is constructed, one can 
choose $\tau$ to be a Nielsen automorphism). 
\end{itemize} 
 
In the present paper we provide a theoretical justification of these 
experimental results. It turns out that the explanation comes from 
exploring the action of $Out(F)$ on the space of geodesic currents on 
$F$, analyzed by the author in~\cite{Ka,Ka1}. 
Recall that $C$ denotes the set of all nontrivial cyclically reduced words in $F$.
 
Our main result is: 
 
\begin{theor}\label{A} 
  Let $F=F(A)$ be a free group where $A=\{a_1,\dots, a_k\}$ and $k\ge 
  2$. Let $\phi\in Aut(F)$ be an arbitrary automorphism that is not a 
  composition of a relabelling and an inner automorphisms.

  Then there exist a Whitehead automorphism $\tau$ of $F$ and a cyclic 
  word $w$ with the following properties: 
 
\begin{enumerate} 
\item For $m_A$-a.e. point $\omega\in \partial F$ we have 
$$ 
||\tau\phi(\omega_n)||<||\phi(\omega_n)|| 
$$ 
as $n\to\infty$ and 
$$ 
\lim_{n\to\infty} [\Gamma_{\phi(\omega_n)}]=[\Gamma_{\phi(w)}], 
$$ 
where $[\Gamma_g]$ is the normalized Whitehead graph corresponding to 
the conjugacy class of $g\in F$. 
 
\item For every $\epsilon>0$ there is a $C$-exponentially generic 
  subset $U\subseteq C$ such that for each $f\in C$ 
$$ 
||\tau\phi(f)||<||\phi(f)|| 
$$ 
and 
$$ 
d\big([\Gamma_{\phi(f)}],[\Gamma_{\phi(w)}]\big)\le \epsilon. 
$$

\item For every $\epsilon>0$ there is an $F$-exponentially generic 
  subset $W\subseteq F$ such that for every $f\in W$ 
$$ 
||\tau\phi(f)||<||\phi(f)|| 
$$ 
and 
$$ 
d\big([\Gamma_{\phi(f)}],[\Gamma_{\phi(w)}]\big)\le \epsilon. 
$$ 
\end{enumerate} 
 
\end{theor} 
 
The definitions of genericity, Whitehead graphs and the uniform 
measure $m_A$ are given in the subsequent sections. Informally, if 
$\omega\in \partial F$ is an $m_A$-random point, the element 
$\omega(n)\in F$ is a ``random'' freely reduced element of length 
$n$, which is also close to being cyclically reduced. Normalized 
Whitehead graphs of a cyclically reduced word $w$, roughly 
speaking, records the frequencies with which the two-letter freely 
reduced words occur in $w$.

Thus Theorem~\ref{A} shows that, in terms of the experiments described 
above, there will be one ``cluster'' for each $\phi\in \Phi$ 
consisting of all $\phi_w(w)$ such that $\phi_w=\phi$, $w\in W_2$.

Note that the Whitehead automorphism $\tau$ in the statement of 
Theorem~\ref{A} is algorithmically computable in terms of $\phi\in 
Aut(F)$, although the complexity of such an algorithm is a priori 
exponential in terms of the word length of the outer automorphism 
$[\phi]$ in $Out(F)$.  The main tool used in the proof of
Theorem~\ref{A} is the machinery of \emph{geodesic currents} on free
groups, discussed in more detail in Section~\ref{sect:curr} below.

Together with a recent result of 
S.~Francaviglia~\cite{Fr}, the proof of Theorem~\ref{A} turns out to imply
the existence of the following ``universal" length-reducing factorization for
automorphisms of free groups (when applied to ``random" elements of $F$):

\begin{corol}\label{B}\cite{Fr} 
  Let $\phi\in Aut(F)$ be an arbitrary automorphism that is not a
  composition of a relabelling and an inner automorphisms.  Then there
  exists a factorization
$$ 
\phi=\sigma_m\sigma_{m-1}\dots \sigma_1 \alpha, 
$$ 
where $m\ge 1$, the automorphism $\alpha$ is a 
  composition of a relabelling and an inner automorphisms, where $\sigma_i$ 
are Whitehead automorphisms of the second kind, and such that the
following holds.

Denote $\psi_0=\alpha$, $\psi_i=\sigma_i\sigma_{i-1}\dots \sigma_1
\alpha$ for $i=1,\dots, m$. Thus $\psi_m=\phi$.

Then for $m_A$-a.e. point $\omega\in \partial F$ as $n\to\infty$ we have
$$
||\psi_i\omega_n||<||\psi_{i+1}\omega_n||, \qquad i=1, \dots, m-1 
$$
so that
$$ 
||\omega_n||=||\psi_0\omega_n||<||\psi_1\omega_n||<\dots <||\psi_m\omega_n||=||\phi\omega_n||.
$$ 

\end{corol}

\section{Geodesic Currents} \label{sect:curr}
 
We recall some basic notions related to geodesic currents on free 
groups. We refer the reader to~\cite{Ka,Ka1,Ma} for a more 
comprehensive discussion.

\begin{conv} 
  We identify the hyperbolic boundary $\partial F$ with the set of all 
  geodesic rays from $1$ in $X$ or equivalently, with the set of all 
  semi-infinite freely reduced words 
$$ 
\omega=a_1a_2\dots a_n\dots, \text{ where } a_i\in A^{\pm 1}. 
$$ 
The boundary $\partial F$ is endowed with the Cantor-set topology and 
with the homeomorphic left $F$-action by left translations, as usual. 
We also denote 
$$ 
\partial^2 F:=\{(\zeta,\xi): \zeta,\xi\in \partial F \text{ and } 
\zeta\ne \xi\}. 
$$ 
Note that $\partial^2 F$ comes equipped with the diagonal left 
$F$-action by homeomorphisms. 
 
For a directed geodesic segment $\gamma=[x,y]$ in $X$ with $x,y\in F$, 
$x\ne y$ we denote by $Cyl_X(\gamma)$ the set of all $(\zeta,\xi)\in 
\partial^2 F$ such that the geodesic $[\zeta,\xi]$ in $X$ passes 
through $\gamma$ in the correct direction. Note that 
$Cyl_X(\gamma)\subseteq\partial^2 F$ is an open-closed compact subset 
of $\partial^2 F$. 
 
We denote by $\mathcal P(X)$ the set of all directed geodesic segments 
of positive length in $X$ with endpoints in $VX=F$. Also, denote 
$F_\ast:=F-\{1\}$. 
\end{conv}

\begin{defn}[Uniform measure] 
  For $v\in F_\ast$ denote by $Cyl_A(v)$ the set of all geodesic rays 
  $\omega\in \partial F$ that begin with $v$. 
 
  The \emph{uniform measure} $m_A$ on $\partial F$ is the Borel 
  probability measure on $\partial F$ defined by 
$$ 
m_A(Cyl_A(v))=\frac{1}{2k(2k-1)^{|v|-1}}\quad \text{ for every } v\in 
F_\ast. 
$$ 
 
\end{defn}

\begin{defn}[Geodesic currents] 
  A \emph{geodesic current} on $F$ is a locally finite (that is finite
  on compact subsets) positive Radon measure $\nu$ on $\partial^2 F$
  such that $\nu$ is $F$-invariant.  The set of all geodesic currents on
  $F$ is denoted by $Curr(F)$.  The space $Curr(F)$ comes equipped with
  the natural weak topology which can be described as follows.
 
  For $\nu_n,\nu\in Curr(F)$ we have 
$$ 
\lim \nu_n =\nu 
$$ 
if and only if 
$$ 
\lim_{n\to\infty}\nu_n(Cyl_X(\gamma))=\nu(Cyl_X(\gamma)) \text{ for 
  every } \gamma\in \mathcal P(X). 
$$ 
\end{defn} 
 
\begin{defn}[The coordinates on $Curr(F)$] 
  If $\nu\in Curr(F)$ and $\gamma=[x,y]\in \mathcal P(X)$ then by 
  $F$-invariance of $\nu$ the value $\nu(Cyl_X(\gamma))$ only depends 
  on $\nu$ and the \emph{label} $v:=x^{-1}y\in F$ of $\gamma$.  For a 
  nontrivial $v\in F$ we denote 
$$ 
\langle v,\nu\rangle:=\nu(Cyl_X(\gamma)) 
$$ 
where $\gamma\in \mathcal P(X)$ is any geodesic segment labelled by 
$v$. 
We call $\langle v,\nu\rangle$ the \emph{number of occurrences of $v$ 
  in $\nu$}. 
\end{defn} 
 
The following lemma~\cite{Ka1} summarizes some basic invariance 
properties satisfied by the coordinates of a geodesic current: 
 
\begin{lem}\label{inv} 
  Let $\nu\in Curr(F)$. Then for every $v\in F_\ast$ 
$$ 
\langle v,\nu\rangle=\sum_{a\in A^{\pm 1}, |va|=|v|+1} \langle 
va,\nu\rangle =\sum_{a\in A^{\pm 1}, |av|=1+|v|} \langle 
av,\nu\rangle. 
$$ 
\end{lem} 
 
A current $\nu\in Curr(F)$ is uniquely determined by the family 
$(\langle v, \nu\rangle)_{v\in F_\ast}$. Moreover, as shown 
in~\cite{Ka,Ka1}, every nonnegative family $(\langle v, 
\nu\rangle)_{v\in F_\ast}$, satisfying the invariance conditions from 
Lemma~\ref{inv}, defines a current $\nu\in Curr(F)$.

\begin{defn}[Uniform current]\label{defn:uc} 
  The \emph{uniform current $n_A\in Curr(F)$} corresponding to the 
  free basis $A$ of $F$ is the geodesic current defined by: 
$$ 
n_A(Cyl_X(\gamma))=\frac{1}{2k(2k-1)^{|\gamma|-1}} \text{ for every } 
\gamma\in \mathcal P(X). 
$$ 
 
Thus $\langle v,n_A\rangle=\frac{1}{2k(2k-1)^{|v|-1}}$ for every $v\in F_\ast$ 
\end{defn}

\begin{defn}[Rational currents] 
  Let $g\in F_\ast$.  If $g$ is not a proper power, define 
$$ 
\eta_g:=\sum_{h\in [g]} \delta_{(h^{-\infty},h^{\infty})}. 
$$ 
where $[g]$ is the conjugacy class of $g$ in $F$. 
If $g=g_0^s$ where $s\ge 2$ and $g_0\in F_\ast$ is not a proper power, 
define 
$$ 
\eta_g:=s \eta_{g_0}. 
$$ 
It is easy to see that $\eta_g$ depends only on the conjugacy class 
$[g]$ of $g$ in $F$. 
 
Nonnegative multiples of the currents $\eta_g, g\in F_\ast$, are 
called \emph{rational currents}. 
\end{defn} 
 
An important basic fact (see~\cite{Ka1}) is: 
\begin{prop} 
  The set of rational currents is dense in $Curr(F)$. 
\end{prop}

\begin{conv}[Cyclic words] 
  We will often think about conjugacy classes of nontrivial elements 
  of $F$ as \emph{cyclic words}. A \emph{cyclic word} $w$ over $A$ is 
  a nontrivial cyclically reduced word in $F(A)$ written clockwise on 
  a circle without specifying an initial point. The length of that 
  cyclically reduced word is called the \emph{cyclic length} of $w$ 
  and is denoted by $||w||$.  The circle is thought of as a labelled 
  graph subdivided into $||w||$ directed edges, each labelled by a 
  letter of $A$. 
 
  If $v\in F$, we call a vertex on this circle an \emph{occurrence of 
    $v$ in $w$} if $v$ can be read in the circle starting at that 
  vertex and going clockwise (we are allowed to stop at a different 
  vertex from the one where we started). The number of occurrences of 
  $v$ in $w$ is denoted by $\langle v, w\rangle$. 
 
  Also, if $v,g\in F$ are nontrivial elements, we put $\langle v, 
  g\rangle:=\langle v, w\rangle$ where $w$ is the cyclic word 
  representing the conjugacy class of $g$. 
\end{conv} 
The following basic fact gives a useful alternative description of 
rational currents:

\begin{lem}\label{lem:rational} 
  Let $g\in F_\ast$ and let $w$ be the cyclic word determined by the 
  conjugacy class of $g$. Then for every $v\in F_\ast$ we have 
$$ 
\langle v,w\rangle=\langle v, \eta_g\rangle. 
$$ 
\end{lem} 
 
There is a natural continuous left action of $Aut(F)$ on $Curr(F)$ 
which factors to the action of $Out(F)$ on $Curr(F)$. If $\phi\in 
Aut(F)$ then $\phi$ is a quasi-isometry of the Cayley graph $X$ of 
$F$. Therefore $\phi$ induces a canonical boundary homeomorphism 
$\partial \phi:\partial F\to\partial F$ which diagonally extends to a 
homeomorphism $\partial^2\phi :\partial^2 F\to\partial^2 F$.  If 
$\nu\in Curr(F)$ and $\phi\in Aut(F)$, the current $\phi\nu\in 
Curr(F)$ is defined by setting 
$$ 
\phi\nu(S):=\nu((\partial^2\phi)^{-1}(S)) 
$$ 
for every Borel subset $S\subseteq \partial^2 F$.  It is not hard to 
show (see~\cite{Ka1}) that for every $g\in F_\ast$ and every $\phi\in 
Aut(F)$ we have $\phi\eta_g=\eta_{\phi(g)}$.

The following useful statement, established in~\cite{Ka1}, gives a 
``coordinate'' description of the action of $Aut(F)$ on $Curr(F)$. 
 
\begin{prop}\label{action} 
  Let $\phi\in Aut(F)$. there exists an integer $K=K(\phi)>0$ with the 
  following property. 
 
  For every $v\in F_\ast$ there exists a collection of nonnegative 
  integers $\{c(u,v,\phi): u\in F, |u|=K|v|\}$ such that for every 
  $\nu\in Curr(F)$ 
$$ 
\langle v,\phi \nu\rangle=\sum_{u\in F, |u|=K|v|} c(u,v,\phi)\langle 
u,\nu\rangle. 
$$ 
 
\end{prop}

If $a_n,a\in \mathbb R$ and $\lim_{n\to\infty} =a$, we say that 
the convergence in this limit is \emph{exponentially fast} if 
there exist $0<\sigma<1$, $b>0$ such that $|a_n-a|\le b \sigma^n$ 
for all $n\ge 1$. 
 
\begin{defn}[Generic sets] 
Let $S\subseteq F$ be an infinite subset. Let $T\subseteq S$. We 
say that $T$ is \emph{generic} in $S$, or \emph{$S$-generic} if 
$$ 
\lim_{n\to\infty} \frac{\#\{g\in T: |g|\le n\}}{\#\{g\in S: |g|\le 
n\}}=1. 
$$ 
If, in addition, the convergence in this limit is exponentially 
fast, we say that $T$ is \emph{exponentially $S$-generic}. 
\end{defn} 
In practice we will only be interested in the cases where $S=F$ or 
$S=C$ (recall that $C$ is the set of all nontrivial cyclically reduced
words in $F$). We refer the reader to~\cite{KMSS,KSS} for more details 
regarding genericity and generic-case complexity.

\section{The length functional} 
 
It turns out that the notion of ``cyclic length'' with respect to the 
free basis $A$ extends naturally to a continuous linear function on  
$Curr(F)$. 
 
\begin{defn}[Length of a current] 
  Let $\nu\in Curr(F)$.  We define the \emph{length $L(\nu)$ of $\nu$ 
    with respect to $A$} as: 
$$ 
L(\nu):=\sum_{a\in A^{\pm 1}} \langle a,\nu\rangle. 
$$ 
\end{defn} 
In the language of~\cite{Ka1} we have $L(\nu)=I(\ell_A,\nu)$ where $I$ 
is the ``intersection form'' and where $\ell_A:F\to\mathbb R$ is the 
length function defined as $\ell_A(w)=||w||$ for $w\in F$.  Note that 
for any automorphism $\phi\in Aut(F)$ the number $L(\phi n_A)$ is 
exactly what in~\cite{KKS} is called the \emph{generic stretching 
  factor} $\lambda_A(\phi)$ of $\phi$ with respect to $A$.

The following basic properties of length follow directly from the 
results about the intersection form established in~\cite{Ka1}.

\begin{prop}\label{prop:length} 
  The following hold: 
\begin{enumerate} 
\item The function $L: Curr(F)\to\mathbb R$ is continuous and linear. 
\item For any integer $m\ge 1$ and for every $\nu\in Curr(F)$ we have 
$$ 
L(\nu)=\sum_{v\in F, |v|=m} \langle v,\nu\rangle. 
$$ 
 
\item For every $w\in F_\ast$ we have 
$$ 
||w||=L(\eta_w). 
$$ 
\item We have $L(n_A)=1$. 
\end{enumerate} 
\end{prop} 
 
In view of Proposition~\ref{action} and Proposition~\ref{prop:length} 
we obtain: 
 
\begin{prop}\label{prop:le} 
  Let $\phi\in Aut(F)$. 
 
\begin{enumerate} 
\item There is $m\ge 2$ and a collection of integers $\{d(u): u\in F, 
  |u|=m\}$ such that for every $\nu\in Curr(F)$ we have 
$$ 
L(\phi\nu)=\sum_{|u|=m} d(u)\langle u,\nu\rangle. 
$$ 
\item Suppose $m\ge 1$ is an integer and $\{d(u)\in \mathbb Z: u\in 
  F,|u|=m\}$ are such that for every cyclic word $w$ we have 
$$ 
||\phi(w)||=\sum_{|u|=m} d(u)\langle u,w\rangle. 
$$ 
Then for every $\nu\in Curr(F)$ we have 
$$ 
L(\phi\nu)=\sum_{|u|=m} d(u)\langle u,\nu\rangle. 
$$ 
\end{enumerate} 
\end{prop} 
 
\begin{proof} 
  Part (1) follows directly from Proposition~\ref{action} and 
  Proposition~\ref{prop:length}. Suppose the assumptions of part (2) 
  hold. Then the conclusion of part (2) holds for every current of the 
  form $\eta_g$, $g\in F_\ast$. Therefore, in view of 
  Lemma~\ref{lem:rational}, the conclusion of part (2) holds for every 
  $\nu\in Curr(F)$ since rational currents are dense in $Curr(F)$. 
\end{proof}

\section{Whitehead automorphisms}

Recall that $\Sigma=A\cup A^{-1}=\{a_1,\dots, a_k, a_1^{-1}, \dots, 
a_k^{-1}\}$.  We follow Lyndon and Schupp, Chapter~I~\cite{LS} in our 
discussion of Whitehead automorphisms.  We recall the basic 
definitions and results.

\begin{defn}[Whitehead automorphisms]\label{defn:moves} 
  A \emph{Whitehead automorphism} of $F$ is an automorphism $\tau$ of 
  $F$ of one of the following two types: 
 
  (1) There is a permutation $t$ of $\Sigma$ such that 
  $\tau|_{\Sigma}=t$. In this case $\tau$ is called a \emph{relabeling 
    automorphism} or a \emph{Whitehead automorphism of the first 
    kind}. 
 
  (2) There is an element $a\in \Sigma$, the \emph{multiplier}, such 
  that for any $x\in \Sigma$ 
$$ 
\tau(x)\in \{x, xa, a^{-1}x, a^{-1}xa\}. 
$$ 
 
In this case we say that $\tau$ is a \emph{Whitehead automorphism of 
  the second kind}. (Note that since $\tau$ is an automorphism of $F$, 
we always have $\tau(a)=a$ in this case). To every such $\tau$ we 
associate a pair $(T,a)$ where $a$ is as above and $T$ consists of all 
those elements of $\Sigma$, including $a$ but excluding $a^{-1}$, such 
that $\tau(x)\in\{xa, a^{-1}xa\}$.  We will say that $(T,a)$ is the 
\emph{characteristic pair} of $\tau$. 
\end{defn} 
 
Note that for any $a\in \Sigma$ the inner automorphism corresponding to
the conjugation by $a$ is a 
Whitehead automorphism of the second kind. 
 
\begin{defn}[Minimal elements] 
  An element $w\in F$ is said to be \emph{automorphically minimal} or 
  just \emph{minimal} if for \emph{every} $\alpha\in Aut(F)$ we have 
  $|w|\le |\alpha(w)|$. 
\end{defn}

\begin{prop}\label{wh}[Whitehead's  Algorithm] 
\begin{enumerate} 
\item If $u \in F$ is cyclically reduced and not minimal, then there 
  is a Whitehead automorphism $\tau$ such that $||\tau(u)||<||u||$. 
\item Let $u,v \in F$ be minimal (and hence cyclically reduced) 
  elements with $|u|=|v|=n>0$. Then $Aut(F)u=Aut(F)v$ if and only if 
  there exists a finite sequence of Whitehead automorphisms 
  $\tau_s,\dots, \tau_1$ such that $\tau_s\dots \tau_1(u)=v$ and such 
  that for each $i=1,\dots, s$ we have 
$$ 
||\tau_i\dots \tau_1(u)||=n. 
$$ 
\end{enumerate} 
\end{prop} 
 
\begin{defn}[Strict Minimality] 
  A nontrivial cyclically reduced word $w$ in $F$ is \emph{strictly 
    minimal} if for every non-inner Whitehead automorphism $\tau$ of 
  $F$ of the second kind we have 
$$ 
||\tau(w)||>||w||. 
$$ 
\end{defn} 
 
\begin{defn}[Simple automorphisms] 
  An automorphism $\phi\in Aut(F)$ is called \emph{simple} if it is 
  the composition of an inner and a relabelling automorphisms. 
\end{defn} 
 
Clearly if $\phi$ is simple, then for every $w\in F_\ast$ we have 
$||\phi(w)||=||w||$. Proposition~\ref{wh} immediately implies that 
every strictly minimal element is minimal and, moreover, if $u$ is 
strictly minimal and $\phi\in Aut(F)$ is such that $||u||=||\phi(u)||$ 
then $\phi$ is simple. 
 
\begin{defn}[Weighted Whitehead graph] 
 
  Let $w$ be a nontrivial cyclic word in $F(A)$.  The \emph{weighted 
    Whitehead graph $\Gamma_w$ of $w$} is defined as follows.  The 
  vertex set of $\Gamma_w$ is $\Sigma$. For every $x,y\in \Sigma$ such 
  that $x\ne y^{-1}$ there is an undirected edge in $\Gamma_w$ from 
  $x^{-1}$ to $y$ labelled by the sum 
$$ 
\langle xy,w\rangle+\langle y^{-1}x^{-1},w\rangle, 
$$ 
the number of occurrences of the words $xy$ and $y^{-1}x^{-1}$ in $w$.

The \emph{normalized Whitehead graph} $[\Gamma_w]$ of $w$ is the 
labelled graph obtained from $\Gamma_w$ by dividing every edge-label 
by $||w||$. 
\end{defn} 
 
\begin{defn} 
  An \emph{abstract Whitehead graph} is a labelled graph $\Gamma$ 
  whose vertex and edge sets are the same as those for a weighted 
  Whitehead graph of a cyclic word and such that each edge $e$ of 
  $\Gamma$ is labelled by a real number $r(e)$. If $\Gamma,\Gamma'$ 
  are two abstract Whitehead graphs, we define 
  $$d(\Gamma,\Gamma')=\max_{e\in E\Gamma} |r(e)-r(e')|.$$ This turns 
  the set of all abstract Whitehead graphs into a metric space 
  homeomorphic to $\mathbb R^{k(2k-1)}$. 
\end{defn} 
Note that if $w$ is a cyclic word, then both $\Gamma_w$ and 
$[\Gamma_w]$ are abstract Whitehead graphs. Note also that for 
$[\Gamma_w]$ the sum of all edge-labels is equal to $1$.

\begin{conv} 
  Let $w$ be a fixed nontrivial cyclic word.  For two subsets 
  $P,Q\subseteq \Sigma$ we denote by $P\underset{w}{.}Q$ the sum of 
  all edge-labels in the weighted Whitehead graph $\Gamma_w$ of $w$ of 
  edges from elements of $P$ to elements of $Q$. Thus for $x\in 
  \Sigma$ the number $x\underset{w}{.}\Sigma$ is equal to the total 
  number of occurrences of $x^{\pm 1}$ in $w$. 
\end{conv}

The next lemma, which is Proposition~4.16 of Ch.~I in \cite{LS}, gives 
an explicit formula for the difference of the lengths of $w$ and 
$\tau(w)$, where $\tau$ is a Whitehead automorphism. 
 
\begin{lem}\label{lem:LS} 
  Let $w$ be a nontrivial cyclically reduced word and let $\tau$ be a 
  Whitehead automorphism of the second kind with the characteristic 
  pair $(T,a)$. Let $T'=\Sigma-T$. Then 
 $$ 
 ||\tau(w)||-||w||=T\underset{w}{.}T'-a\underset{w}{.}\Sigma. 
 $$ 
\end{lem}

Now Lemma~\ref{lem:LS} and Proposition~\ref{prop:le} immediately 
imply: 
 
\begin{cor}\label{cor:wh} 
  Let $\tau$ be a Whitehead automorphism of the second kind. Then 
  there exists a collection of integers $\{b(z): z\in F, |z|=2\}$ such 
  that for every $\nu\in Curr(F)$ we have 
$$ 
L(\tau\nu)=\sum_{|z|=2} b(z)\langle z,\nu\rangle. 
$$ 
\end{cor} 
 
\begin{rem} 
Note that in view of Lemma~\ref{lem:LS} and Corollary~\ref{cor:wh}, if 
$w$ is a cyclic word and $\tau$ is a Whitehead automorphism of the 
second kind, then the quantity 
 
$$ 
\frac{||\tau(w)||}{||w||} 
$$ 
is completely determined by $\tau$ and the normalized Whitehead graph 
$[\Gamma_w]$ of $w$. 
\end{rem}

\section{Proof of the main result}

\begin{prop}\label{euler} 
  For every integer $m\ge 2$ there exists a cyclic word $w$ such that 
$$ 
\langle v,w\rangle =1 \text{ for every } v\in F \text{ with } |v|=m 
$$ 
and that $||w||=2k(2k-1)^{m-1}$. 
\end{prop} 
 
\begin{proof} 
 
  This follows from a more general result in \cite{Ka}. We present an 
  argument here for completeness. 
 
  If $v\in F$ is a freely reduced word with $|v|\ge 2$, we denote by 
  $v_{-}$ the initial segment of $v$ of length $|v|-1$ and we denote 
  by $v_+$ the terminal segment of $v$ of length $|v|-1$.

  Let $n\ge 2$.  Form a finite directed labelled graph $\Gamma$ as 
  follows. The vertex set of $\Gamma$ is 
$$ 
V\Gamma:=\{u\in F: |u|=m-1\}. 
$$ 
The set of directed edges of $\Gamma$ is 
$$ 
E\Gamma:=\{v\in F: |v|=m\}. 
$$ 
For each $v\in E\Gamma$ the initial vertex of $v$ in $\Gamma$ is 
$v_{-}$ and the terminal vertex of $v$ in $\Gamma$ is $v_+$. Also, the 
edge $v\in E\Gamma$ is labelled by the \emph{label} $a(v)\in A^{\pm 
  1}$ which is the last letter of the word $v$.

Note that for every vertex $u\in V\Gamma$ both the out-degree of $u$ 
and the in-degree of $u$ in $\Gamma$ are equal to $2k-1$. Thus 
$\Gamma$ is a strongly connected directed graph where for each vertex 
the in-degree is equal to the out-degree. Therefore there exists an 
Euler circuit $c$ is $\Gamma$, that is, a cyclic path passing through 
each directed edge of $\Gamma$ exactly once.  Let $c$ be represented 
by the edge-path 
$$ 
v_1v_2\dots v_{t}, \text{ where } t=|E\Gamma|=2k(2k-1)^{m-1}. 
$$ 
Let $w$ be the cyclic word defined by the word 
$$ 
a(v_1)a(v_2)\dots a(v_t). 
$$ 
 
Then it is not hard to see that $||w||=t=2k(2k-1)^{m-1}$ and that for 
every $v\in F$ with $|v|=m$ we have 
$$ 
\langle v,w\rangle =1, 
$$ 
as required. 
\end{proof} 

Recall that $n_A$ is the uniform current on $F$ defined in Definition~\ref{defn:uc}.

\begin{prop}\label{prop:wm}[Ideal Whitehead Algorithm] 
  Let $\phi\in Aut(F)$ be an automorphism such that $\phi$ is not 
  simple. Then there exists a Whitehead automorphism $\tau$ of the 
  second kind such that 
$$ 
1=L(n_A)\le L(\tau\phi n_A)<L(\phi n_A). 
$$ 
\end{prop}

\begin{proof} 
 
  By~Proposition~\ref{action} there exist an integer $m\ge 2$ and a 
  collection of nonnegative integers 
$$ 
\{c(v,z): v,z\in F, |v|=m, |z|=2\} 
$$ 
such that for every $\nu\in Curr(F)$ we have 
$$ 
\langle z,\phi\nu\rangle =\sum_{|v|=m}c(v,z) \langle v,\nu\rangle. 
$$ 
Let $w$ be a cyclic word provided by Proposition~\ref{euler}. 
Recall that we have $||w||=2k(2k-1)^{m-1}$.  Let 
$\theta=\frac{\eta_w}{2k(2k-1)^{m-1}}$. Thus for every $v\in F$ 
with $|v|=m$ we have 
$$ 
\langle v,\theta\rangle=\frac{1}{2k(2k-1)^{m-1}}. 
$$ 
 
Then for every $z\in F$ with $|z|=2$ we have 
$$ 
\langle z,\phi\theta\rangle=\langle z,\phi n_A\rangle. 
$$ 
 
Moreover, we have 
 
$$ 
L(\phi\theta)=\sum_{|z|=2}\langle z,\phi\theta\rangle= 
\sum_{|z|=2}\langle z,\phi n_A\rangle=L(\phi n_A). 
$$

By Lemma~4.8 of~\cite{KSS} the word $w$ is strictly minimal which 
implies, in particular, that $||w||<||\phi(w)||$, since $\phi$ is not 
simple.  Therefore, by Whitehead's theorem, part (1) of 
Proposition~\ref{wh}, there exists a Whitehead automorphism $\tau$ of 
the second kind such that 
$$ 
||w||\le ||\tau\phi(w)||<||\phi(w)||. 
$$ 
Therefore 
$$ 
1=L(\theta)\le L(\tau \phi\theta)< L(\phi\theta). 
$$ 
then by the above formulas and Corollary~\ref{cor:wh} we see that 
$$ 
1=L(n_A)\le L(\tau\phi n_A)=L(\tau \theta)<L(\theta)=L(\phi n_A), 
$$ 
as required. 
\end{proof} 
 
Note that Proposition~\ref{prop:wm} means that $n_A\in Curr(F)$ is 
``minimal" and even ``strictly minimal" in the sense that for 
every $\phi\in Aut(F)$ 
$$ 
L(n_A)\le L(\phi n_A) 
$$ 
with the equality achieved if and only if $\phi$ is simple.

\begin{cor}\label{cor:ideal} 
Let $\phi\in Aut(F)$ be a non-simple automorphism. Then there 
exists a factorization 
$$ 
\phi=\sigma_m\sigma_{m-1}\dots \sigma_1 \alpha, 
$$ 
where $m\ge 1$, the automorphism $\alpha$ is simple, $\sigma_i$ 
are Whitehead automorphisms of the second kind and 
$$ 
L(\sigma_{i-1}\dots \sigma_1\alpha n_A)< L(\sigma_i\sigma_{i-1}\dots \sigma_1
\alpha n_A), \qquad i=1,\dots, m-1. 
$$ 
\end{cor} 
\begin{proof} 
 
Put 
$$ 
\Lambda:=\{L(\psi n_A): \psi\in Aut(F)\}. 
$$ 
A recent theorem of S.~Francaviglia~\cite{Fr} shows that $\Lambda$ is 
a discrete subset of $\mathbb R$. Also, as proved in~\cite{KSS}, 
for every $\psi\in Aut(F)$ we have $L(\psi n_A)\ge 1$ and, 
moreover, $L(\psi n_A)= 1$ if and only if $\psi$ is simple. 
 
Let $\phi\in Aut(F)$ be a non-simple automorphism. Thus $L(\phi 
n_A)>1$. Repeatedly applying Proposition~\ref{prop:wm} we conclude that 
there exists a sequence of Whitehead automorphisms 
$\tau_1,\tau_2,\dots $ such that $L(\phi n_A)> L(\tau_1\phi 
n_A)>L(\tau_2\tau_1\phi n_A)>\dots$. Since $\Lambda$ is a discrete 
subset of $[1,\infty)$, the sequence $\tau_1,\tau_2,\dots $ must 
terminate in a finite number of steps with some $\tau_m$. Hence 
the automorphism $\alpha:=\tau_m\dots \tau_2\tau_1\phi$ must be 
simple since otherwise by Proposition~\ref{prop:wm} the sequence of $\tau_i$ 
could be extended. Then the factorization 
$$ 
\phi=\tau_1^{-1}\dots \tau_m^{-1}\alpha 
$$ 
has the required properties and the corollary is proved. 
\end{proof}

\begin{proof}[Proof of Theorem~\ref{A}] 
 
  Let $\phi\in Aut(F)$ be an automorphism such that $\phi$ is not simple. By 
  Proposition~\ref{prop:wm} there exists a Whitehead automorphism 
  $\tau$ such that 
 
$$ 
L(\tau\phi n_A)<L(\phi n_A). 
$$ 
 
Also, as in the proof of Proposition~\ref{prop:wm}, let $w$ be the 
cyclic word provided by Proposition~\ref{euler}.

Recall that by Proposition~\ref{action} there exist an integer $m\ge 
2$ and a collection of nonnegative integers 
$$ 
\{c(v,z): v,z\in F, |v|=m, |z|=2\} 
$$ 
such that for every $\nu\in Curr(F)$ we have 
$$ 
\langle z,\phi\nu\rangle =\sum_{|v|=m}c(v,z) \langle v,\nu\rangle. 
$$

Let $\omega\in \partial F$ be an $m_A$-random point.  Then, as 
observed in~\cite{Ka1} 
$$ 
\lim_{n\to\infty} \frac{\eta_{\omega_n}}{n}=\lim_{n\to\infty} 
\frac{\eta_{\omega_n}}{||\omega_n||}=n_A \text{ in } Curr(F). 
$$ 
 
Hence 
$$ 
\lim_{n\to\infty} \phi \frac{\eta_{\omega_n}}{n}=\phi n_A 
$$ 
and 
$$ 
\lim_{n\to\infty} \tau\phi \frac{\eta_{\omega_n}}{n}=\tau\phi n_A 
$$

Since $L:Curr(F)\to \mathbb R$ is continuous, and $L(\tau\phi 
n_A)<L(\phi n_A)$, it follows that for $n\to\infty$ 
 
$$ 
L(\tau\phi \frac{\eta_{\omega_n}}{n})<L(\phi 
\frac{\eta_{\omega_n}}{n}), 
$$ 
Then for $n\to\infty$ 
$$ 
\frac{||\tau\phi(\omega_n)||}{n}<\frac{||\phi(\omega_n)||}{n} 
$$ 
and therefore 
$$ 
||\tau\phi(\omega_n)||<||\phi(\omega_n)||, 
$$ 
as required. 
 
We have seen in (\ddag) that 
$$ 
\langle z,\phi\theta\rangle=\langle z,\phi n_A\rangle \text{ for each 
} z\in F \text{ with } |z|=2 
$$ 
where $\theta=\frac{\eta_w}{2k(2k-1)^{m-1}}$ and 
$||w||=2k(2k-1)^{m-1}$.  Since \[\lim_{n\to\infty} \phi 
\frac{\eta_{\omega_n}}{n}=\phi n_A,\] this implies that for each $z\in 
F$ with $|z|=2$ we have 
$$ 
\lim_{n\to\infty} \frac{\langle z, \phi(\omega_n)\rangle}{n}=\langle 
z,\phi n_A\rangle=\langle z,\phi\theta\rangle=\frac{\langle 
  z,\phi(w)\rangle}{2k(2k-1)^{m-1}}. 
$$

We also have 
$$ 
||\phi(w)||=\sum_{|z|=2} \langle z,\phi(w)\rangle 
=\sum_{|z|=2}\sum_{|v|=m} c(v,z)\langle v, 
w\rangle=\sum_{|z|=2}\sum_{|v|=m} c(v,z) 
$$ 
and 
$$ 
||\phi(\omega_n)||=\sum_{|z|=2} \langle z,\phi(\omega_n)\rangle 
=\sum_{|z|=2}\sum_{|v|=m} c(v,z)\langle v, \omega_n\rangle. 
$$ 
Therefore for any $z'\in F$ with $|z'|=2$ we have 
$$ 
\frac{\langle z', \phi(w)\rangle}{||\phi(w)||}=\frac{\sum_{|v|=m} 
  c(v,z')}{\sum_{|z|=2}\sum_{|v|=m} c(v,z)} 
$$ 
and 
$$ 
\frac{\langle z', 
  \phi(\omega_n)\rangle}{||\phi(\omega_n)||}=\frac{\sum_{|v|=m} 
  c(v,z')\langle v, \omega_n\rangle}{\sum_{|z|=2}\sum_{|v|=m} 
  c(v,z)\langle v, \omega_n\rangle}=\frac{\sum_{|v|=m} 
  c(v,z')\frac{\langle v, 
    \omega_n\rangle}{n}}{\sum_{|z|=2}\sum_{|v|=m} c(v,z)\frac{\langle 
    v, \omega_n\rangle}{n}}. 
$$ 
 
Since $\lim_{n\to\infty} \frac{\eta_{\omega_n}}{n}=n_A$, it follows 
that $\lim_{n\to\infty}\frac{\langle v, 
  \omega_n\rangle}{n}=\frac{1}{2k(2k-1)^{m-1}}$ for every $v\in F$ 
with $|v|=m$.  Therefore for every $z'\in F$ with $|z'|=2$ we have 
$$ 
\lim_{n\to\infty} \frac{\langle z', 
  \phi(\omega_n)\rangle}{||\phi(\omega_n)||}=\frac{\sum_{|v|=m} 
  c(v,z')}{\sum_{|z|=2}\sum_{|v|=m} c(v,z)}=\frac{\langle z', 
  \phi(w)\rangle}{||\phi(w)||}. 
$$ 
It follows that $\lim_{n\to\infty} 
[\Gamma_{\phi(\omega_n)}]=[\Gamma_{\phi(w)}]$ as required.

This establishes part (1) of Theorem~\ref{A}.

Recall that by Proposition~6.2 of~\cite{KSS} if $U\subseteq C$ is an 
exponentially $C$-generic subset, then the set $W$ consisting of all 
$w\in F$ whose cyclically reduced forms are in $U$, is exponentially 
$F$-generic. Therefore part (2) of Theorem~\ref{A} implies part (3). 
 
Thus it remains to prove part (2) of Theorem~\ref{A}.

For any $\epsilon'>0$ define 
$$ 
U(\epsilon')=\{u\in C: \left|\frac{\langle 
  v,u\rangle}{||u||}-\frac{1}{2k(2k-1)^{m-1}}\right|\le \epsilon' 
\text{ for every } v\in F, |v|=m\}. 
$$

Recall also that there exists a collection of integers $\{d(z): z\in 
F, |z|=2\}$ such that 
$$ 
L(\tau\nu)-L(\nu)=\sum_{|z|=2}d(z)\langle z,\nu\rangle\quad \text{ for 
  every }\nu\in Curr(F). 
$$ 
Since $L(\tau\phi n_A)<L(\phi n_A)$, there is $\epsilon''>0$ such that 
for every $\nu\in Curr(F)$ satisfying 
$$ 
|\langle z,\nu\rangle-\langle z,\phi n_A\rangle|\le \epsilon'' 
$$ 
for every $z\in F$ with $|z|=2$ we have $L(\tau\nu)-L(\nu)<0$.

The properties of $c(v,z)$ listed above imply that there is 
$\epsilon'>0$ such that for every $u\in U(\epsilon')$ and for every 
$z\in F,|z|=2$ we have 
$$ 
\left|\langle z,\phi\frac{\eta_u}{||u||}\rangle-\langle z,\phi 
n_A\rangle\right|\le \epsilon''. 
$$ 
Hence for every $u\in U(\epsilon')$ 
$$ 
L(\tau\phi\frac{\eta_u}{||u||})-L(\phi\frac{\eta_u}{||u||})<0, 
$$ 
that is 
$$ 
\frac{||\tau\phi(u)||}{||u||}<\frac{||\phi(u)||}{||u||}\quad\Rightarrow\quad ||\tau\phi(u)||<||\phi(u)||. 
$$ 
 
The set $U(\epsilon')\subseteq C$ is exponentially $C$-generic, as was 
observed in~\cite{KKS}.  The proof of the Whitehead graph assertion of 
part (2) of Theorem~\ref{A} is similar to that used in part (1). One 
shows that if $\epsilon>0$ is arbitrary then for $\epsilon'>0$ small 
enough 
$$ 
d([\Gamma_{\phi(u)}],[\Gamma_{\phi{w}}])\le \epsilon \text{ for all } 
u\in U(\epsilon'). 
$$ 
We leave the details to the reader. 
 
This completes the proof of Theorem~\ref{A}.

\end{proof}

Corollary~\ref{cor:ideal} implies Corollary~\ref{B} from the
Introduction in a way similar to the proof of part (1) of
Theorem~\ref{A} and we leave the details to the reader.
 
\section{Acknowledgements} 
 
The author is grateful to Alexey G. Myasnikov for suggesting to 
consider the question addressed by Theorem~\ref{A}. The author 
especially thanks Paul Schupp for useful conversations and for 
help with computer experimentation. The author was supported by 
the NSF grant DMS-0404991.

Department of Mathematics, University of Illinois at 
Urbana-Champaign, 1409 West Green Street, Urbana, IL 61801, USA; 
kapovich@math.uiuc.edu

\end{document}